\documentclass[10pt]{article}
\usepackage[utf8]{inputenc}
\usepackage[english]{babel}
\usepackage{amsmath,amsthm,amsfonts}
\usepackage{graphicx}
\usepackage{bm}
\usepackage{tikz}
\usepackage{multirow}

\renewcommand{\div}{\mathop{\rm div}\nolimits}

\newcommand{\norm}[1]{\left\|#1\right\|}

\newcommand{\seminorm}[1]{\left|#1\right|}

\begin{document}

\title{
A Generalized Multiscale Finite Element Method for Poroelasticity Problems II: Nonlinear Coupling
\thanks{This work was supported by RFBR (project N 15-31-20856)}} 

\author{
Donald L. Brown
\thanks{
School of Mathematical Sciences,
The University of Nottingham,
University Park
Nottingham, NG7 2RD, United Kingdom
Email: {\tt donaldbrowdr@gmail.com}.
}
\and
Maria Vasilyeva
\thanks{
Department of Computational Technologies, Institute of Mathematics and Informatics, North-Eastern Federal University, Yakutsk, Republic of Sakha (Yakutia), Russia, 677980 \& Institute for Scientific Computation, Texas A\&M University, College Station, TX 77843. 
Email: {\tt vasilyevadotmdotv@gmail.com}.
}
}

\maketitle

\begin{abstract}
In this paper, we consider the numerical solution of some nonlinear poroelasticity problems that 
are of Biot type and develop a general algorithm for solving nonlinear coupled systems. 
We discuss the difficulties associated with flow and mechanics in heterogenous media with nonlinear coupling. 
The central issue being how to handle the nonlinearities and the multiscale scale nature of the media.
 To compute an efficient numerical solution we develop and
implement a  Generalized Multiscale Finite Element Method (GMsFEM)  that solves nonlinear problems on a
coarse grid by constructing local multiscale basis functions and treating part of the nonlinearity locally as a parametric value.
 After linearization with a Picard Iteration,
 the procedure begins with construction
of multiscale bases for both displacement and pressure in each coarse block by treating the staggered nonlinearity as a parametric value. 
 Using a snapshot space 
and local spectral problems, we construct an offline basis of reduced dimension. 
 From here an online, parametric dependent, space is constructed.
 Finally, after multiplying by
a multiscale partitions of unity, the multiscale basis is constructed and the coarse
grid problem then can be solved for arbitrary forcing and boundary conditions.
  We implement this
algorithm on a geometry with a linear and nonlinear pressure dependent permeability field and compute error between the multiscale solution with the fine-scale solutions. 
\end{abstract}


\section{Introduction}

The applications of mechanics and flow in porous media are wide ranging, as are the challenges involved in simulating some of these problems in nonlinear multiscale contexts.
This is particularly true in geomechanical modeling where relevant phenomena may be highly nonlinear, for example in the setting of 
 enhanced production and environmental safety concerns due to overburden subsidence and compaction \cite{sayers2007introduction,zoback2010reservoir}. 
Another of the central challenges is the multiscale nature of the media considered in geomechanics problems.
 Heterogeneity of rock properties should be accurately accounted in the geomechanical model, and this  requires a computationally costly a high resolution solve.
  Moreover, due to the multi-physics nature of the problems, they may involve highly nonlinear relations. This then makes the further requirement of many iterations in a Newton or Picard  linearization. 
 Thus, we propose a multiscale method  to attempt overcome some of these challenges. The central idea is to linearize in a Picard iteration, and treat the nonlinearities as a parametric value as utilized in \cite{nonlinGMsFEM} and references therein. 
%

As noted in \cite{BrownVasilyeva}, the basic mathematical structure  of the poroelasticity models  are usually coupled equations for pressure and displacements known as Biot type models \cite{BiotOriginal}. 
The pressure equations are a parabolic equation coupled to a time derivative of volumetric strain. While the mechanics equations are are given by quasi-static elasticity equations and is coupled by gradients of pressure. 
In this work however, we focus on the possible nonlinear couplings of the Biot model. There are a myriad of physical and modeling  reasons  to add nonlinearity to the Biot equations, however, we will primarily focus when the permeability field and elasticity tensors depend nonlinearly on pressure and displacements and their gradients. This is due primarily to us wanting to focus on the nonlinearities effects on our GMsFEM, as  nonlinearities in lower order derivative will not interfere with the construction of the local multiscale basis functions.

Nonlinear Poroelastic models of this type have been explored in the  literature to incorporate higher order physics considerations. For example, when the viscosity of the fluid heavily depends on the fluid pressure we may obtain relations of permeability of the form 
$$K(x;p)=\frac{k(x)}{\mu(x;p)}.$$
Here $k$ is the absolute permeability, and $\mu(x;p)$ is the pressure and spatially  dependent viscosity. This can occur when their are very high pressure gradients \cite{Rajagopal}. In the setting of complex geomechanical interactions \cite{minkoff2003coupled} used a relationship between permeability and volumetric strain of the form
$$
K(x;\nabla \cdot u)=A(x)\exp(B(x)\nabla \cdot u),
$$
where $A,B$ are determined constants and $\nabla \cdot u$ is the volumetric strain. Further in \cite{minkoff2003coupled}, the porosity $\phi$  also depends linearly on $(\nabla p, \nabla \cdot u),$ however, this is multiplied throughout generating a nonlinearity.
 In the context of fractured reservoirs, permeability is often computed via the so called "cubic-law" through channels and this may be coupled in orientation and magnitude via the displacements in a nonlinear way \cite{settarifractured}. With this GMsFEM, we propose a method to efficiently compute solutions to these nonlinear poroelasticity problems  with the heterogeneous multiscale properties.

As noted in the prequel \cite{BrownVasilyeva}, there are many very effective multiscale frameworks that have been developed in recent years. 
There are rigorous approaches based on homogenization of partial differential equations \cite{brown2014,brown2011}. However, these approaches may have limited computational use. Examples of computational approaches  include the Heterogeneous Multiscale Method (HMM) \cite{E:Engquist:2003,Abdulle:E:Engquist:Vanden-Eijnden:2012},
  approaches based on the Variational Multiscale Method (see \cite{MR2300286}), where coarse-grid quasi-interpolation operators are used to build an orthogonal splitting into a multiscale space and a fine-scale space \cite{MP11}. 
In this paper, 
we will use the Generalized Multiscale Finite Element Method and its extension to nonlinear poroelasticity  problems in the framework of  \cite{nonlinGMsFEM}.
Specifically, to handle multiscale nonlinear problems, we combine ideas of model reduction, whereby the nonlinearity is replaced locally by a parameter space and offline and online spaces are generated. For a broad presentation of these methods we refer the reader to \cite{egh12}.

%


The paper is organized similarly to \cite{BrownVasilyeva},  as follows. In Section 2 we provide the mathematical background of the nonlinear poroelasticity problem. We will introduce the Biot type model and highlight where the heterogeneities primarily occur. We again note that the nonlinearities in our model are in the permeability and elasticity tensor as these are second order derivative terms.
In  Section 3, to outline the difficulties in full direct numerical simulation we introduce the fine-scale discretizations using coupled time-stepping schemes and a Picard iteration technique for linearization.   In Section 4, we present our nonlinear GMsFEM algorithm and outline its construction procedure. 
Finally, numerical implementations are presented in Section 5. Using the GMsFEM,  we compare the multiscale solution to fine-scale solutions and give error estimates. 
We will present two different examples with permeability being linear and nonlinear with respect to pressure. 
Additionally, we will implement and discuss different snapshot spaces and coarse-grids choices, and its relation to enrichment and the error. 

\section{Problem formulation} 
We denote our computational domain $\Omega\subset\mathbb{R}^d$ to be a bounded Lipschitz region.  We 
 consider a general nonlinear poroelastic system where we wish to find a pressure $p$ and displacements $u$ satisfying 
 \begin{subequations}\label{eq:main}
\begin{eqnarray}
\label{eq:poroelas}
- \div  \left(C (x; u,\nabla u,p,\nabla p)\varepsilon(u)\right) - \alpha \nabla  p  &= 0 \text{ in } \Omega, \\
- \div \left( K(x; u,\nabla u,p,\nabla p) \nabla p \right) + \alpha \frac{\partial \div  u}{\partial t} + \frac{1}{M} \frac{\partial p}{\partial t} &= f \text{ in } \Omega,
\end{eqnarray}
\end{subequations}
with initial condition for pressure 
$
p( x, 0) = p_0. 
$
We write the boundary of the domain into four sections  $\partial \Omega = \Gamma_1 \cup \Gamma_2 = \Gamma_3  \cup \Gamma_4$. We suppose the following boundary conditions on each portion
\[
 C (x; u,\nabla u,p,\nabla p)\varepsilon(u)\cdot  n = 0, \quad  x \in \Gamma_1, \quad 
 u = u_1, \quad  x \in \Gamma_2, 
 \]
 and 
 \[
-K (x; u,\nabla u,p,\nabla p) \frac{\partial p}{\partial n}= 0, \quad  x \in \Gamma_3, \quad  
p = p_1, \quad  x \in \Gamma_4.
\]
Here the symmetric strain is written as  $\varepsilon(u) = \frac{1}{2} \left( \nabla  u + \nabla  u^T \right)$ and we write 
$$C (x; u,\nabla u,p,\nabla p)\varepsilon(u):=C (x; u,\nabla u,p,\nabla p):\varepsilon(u)$$ to mean the double contraction of a $4$-tensor with a $2$-tensor. 

As in the linear case, the primary sources of the heterogeneities in the physical properties arise from $ C (x; u,\nabla u,p,\nabla p)$, the elastic tensor, and $K (x; u,\nabla u,p,\nabla p)$, the permeability. In this setting, we suppose these heterogenous parameters can depend in $p$ and $u$ and their gradients in complicated nonlinear ways. 
Further, we will denote $M$ to be the Biot modulus, $\nu$ is the fluid viscosity, and $\alpha$ is the Biot-Willis fluid-solid coupling coefficient.  Here,  $f$ is a source term representing injection or production processes and $n$ is the unit normal to the boundary. Body forces, such as gravity, are neglected.
 \vspace{.1in}
 
 {\bf Remark}:
 Note that one could also add nonlinearities in the coefficients $\alpha$ and $M$, however, these correspond to lower order terms with respect to derivatives. Therefore, these will not contribute to the local problems in the GMsFEM.  Hence, we will consider them to be constant throughout.
 
  \vspace{.1in}
 
We recall the setting when these relations become linear. In the case of a linear elastic stress-strain constitutive relation we have that the stress tensor and symmetric strain gradient  may be expressed as
\[
C \varepsilon(u)= 2 \mu  \varepsilon( u) + \lambda \div(u) \,  \mathcal{I},
\]
where $\mu$, $\lambda$ are Lame coefficients, $\mathcal{I}$ is the identity tensor.
Note here this $\mu$ is not to be confused with what is often used as a parameter.
Above we have absorbed into the nonlinear permeability coefficient the fluid viscosity $\nu$, and in the case of linear permeability, we have $$K=\frac{k}{\nu},$$ $k$ being absolute permeability.

The nonlinear poroelasticity problem \eqref{eq:main}, 
 can be written in operator matrix form:
\begin{eqnarray}
\label{eq:mmd}
A ( u,p) + \alpha G p = 0,  \\
\frac{d}{dt} \left( S \, p + \alpha D  u \right) + B(u, p) = f,
\end{eqnarray}
where
\[
A ( u,p) = - \div  \left(C (x; u,\nabla u,p,\nabla p)\varepsilon(u)\right),
\quad 
B(u, p) = - \div \left( K(x; u,\nabla u,p,\nabla p) \nabla p \right),
\]
and $G$ and $D$ are gradient and divergence operators and $S=\frac{1}{M}{\cal I}$.

\section{Fine-Scale Discretization}
We will now present fine-scale approximation and nonlinear solution methods  for the above system. 
%
We will motivate the need for a  multiscale method due to the nonlinearity and the heterogeneity of the poroelasticty problem. 
To approximate  the   solution to \eqref{eq:main} on fine-scale grid we will utilize a standard finite element method. 
The corresponding nonlinear variational form of the continuous problem written as 
\begin{eqnarray}
\label{eq:canu}
a( u, p,  v) + g(p,  v)&= 0, \quad \text{for all}   \, \,  v \in \hat{ V}, \\
d \left( \frac{d  u}{dt}, q \right)  + c \left( \frac{d p}{dt}, q \right) 
+ b(u, p , q)&= (f , q), \quad \text{for all}  \, \,  q \in \hat Q.
\end{eqnarray}
for $u \in V$, $p \in Q$ where 
\[
V  = \{ v \in [H^1(\Omega)]^d:  v( x) = u_1,  x \in \Gamma_2  \}, 
\quad
Q  = \{  q \in H^1(\Omega):  q(x) = p_1,  x \in \Gamma_4 \},
\]
and the test spaces  with  homogeneous boundary conditions are given by
\[
\hat{V} = \{ v \in [H^1(\Omega)]^d:  v(x) = 0,  x \in \Gamma_2  \}, 
\quad
\hat{Q} = \{ q \in H^1(\Omega): q(x) = 0,  x \in \Gamma_4 \}.
\]
We define the following nonlinear forms 
\begin{subequations}\label{nonlinearforms}
\begin{align}
a( u, p,  v) & = \int_{\Omega} (C (x; u,\nabla u,p,\nabla p)\varepsilon(u) , \varepsilon(v) ) dx,\\
b(u, p, q)& = \int_{\Omega} \left( K (x; u,\nabla u,p,\nabla p)\nabla p, \nabla q \right) dx,  
\end{align}
\end{subequations}
and  bilinear and linear forms
\[
c(p, q) = \int_{\Omega} \frac{1}{M} \, p \, q \, dx, 
\quad
 g(p,  v)  = \int_{\Omega}\alpha (\nabla p,  v) dx,
\]
and 
\[
d( u, q) = \int_{\Omega} \alpha \div  u \, q \, dx , 
\quad
(f, q) = \int_{\Omega} f  \, q \, dx.
\]
Here $\left( \cdot, \cdot \right)$ under the integrand denotes the standard inner product. In Section \ref{Numerics}, we will discretize the spaces using a fine-scale standard FEM  and denote them $V_{h}, Q_{h}$ and $\hat{V}_{h}, \hat{Q}_{h}$, $h$ being  the fine-grid size. 
The FEM using these spaces will serve as a reference solution for our GMsFEM outlined in Section \ref{GMsFEM}.

{\bf {  Nonlinear Solve:}} We will first consider the time discretizations of the above system, then will discuss resolving the nonlinearity.
 This discretization leads to several possible couplings between time-steps and the two equations of linear poroelasticity \cite{kim2011,vab2014}. However in the nonlinear case we will only consider the fully coupled scheme.
 We proceed by introducing for the nonlinear fully coupled  time derivative operators and then the Picard iteration for the linearization of the nonlinear operators. 

The standard fully implicit finite-difference scheme, or coupled scheme, can be used for the time-discretization and is given by 
\begin{subequations}\label{coupled}
\begin{eqnarray}
\label{eq:capp}
a( u^{n+1},p^{n+1},  v) + g(p^{n+1},  v) = 0,  \\
d \left( \frac{ u^{n+1} -  u^n}{\tau}, q \right)  
+ c \left( \frac{p^{n+1} - p^n}{\tau}, q \right)
+ b(u^{n+1}, p^{n+1} , q) = (f , q),  
\end{eqnarray}
\end{subequations}
with $u^n =  u(x, t_n)$,  $p^n = p(x, t_n)$, where $t_n = n \tau$, $n = 0, 1, ..., M_{T}$, $M_{T} \tau = T$ and $\tau > 0$.
%


We will now consider nonlinear solve in space after time discretization by  the fully coupled scheme \eqref{coupled}. 
%
One could rewrite \eqref{coupled} as a nonlinear system each time step and use a Newton solver, however, for our GMsFEM we prefer to use a linearization based on Picard iteration. Indeed, we may linearize \eqref{nonlinearforms} given $(p_j, u_j)$ from a previous iteration step we write 
\begin{align*}
a( u_{j+1}, p_{j+1},  v)  &\approx a^L( \mu_j ;u_{j+1},  v) :=  \int_{\Omega} (C (x; u_{j},\nabla u_{j},p_{j},\nabla p_{j})\varepsilon(u_{j+1}) , v)dx, 
\\
b(u_{j+1}, p_{j+1}, q) &\approx b^L(\mu_j ; p_{j+1},  v) := \int_{\Omega} \left( K (x; u_{j},\nabla u_{j},p_{j},\nabla p_{j})\nabla p_{j+1}, \nabla q \right) dx,
\end{align*}
where $\mu_j=(u_{j},\nabla u_{j},p_{j},\nabla p_{j})$. We choose this notation in part to emphasize this may be viewed as a parameter in offline phase of the GMsFEM.

Fixing the time-step at $(n+1)$, and 
taking $\mu_{j}=({u}^{n+1}_{j},\nabla {u}^{n+1}_{j},{p}^{n+1}_{j},\nabla {p}^{n+1}_{j})$,  as data from the previous iteration. For $j=0,1,2,\dots$, we wish to find $(u^{n+1}_{j+1},p^{n+1}_{j+1})$ such that 
\begin{subequations}\label{coupled.linear}
\begin{eqnarray}
a^L( \mu_j ;u^{n+1}_{j+1},  v) + g(p^{n+1}_{j+1},  v) = 0,  \\
d \left( \frac{ u_{j+1}^{n+1} -  u^n}{\tau}, q \right)  
+ c \left( \frac{p^{n+1}_{j+1} - p^n}{\tau}, q \right)
+  b^L(\mu_j ; p^{n+1}_{j+1},  q)  = (f , q),  
\end{eqnarray}
\end{subequations}
Once  the desired convergence criteria is reached, we can set the terminal  $(u^{n+1}_{j},p^{n+1}_{j})$ as previous time data.  We  then return to the algorithm time-stepping and continue the iterative  linearization until the terminal time. Note that this process can also be used in an appropriate nonlinear generalization to a  fixed stress splitting \cite{kim2011,vab2014}.

\section{GMsFEM for nonlinear poroelasticity problem}\label{GMsFEM}
We will present the offline and online multiscale basis construction in the fluid or pressure solve then its construction in the mechanics or displacement calculation step in this nonlinearly coupled formulation. Similar to the presentation outlined in \cite{BrownVasilyeva}, however, we will focus on the effects of the nonlinearities on the method. Observing   the linearized formulation \eqref{coupled.linear}, we see  that we may consider the nonlinearity as parametric values we are able to successful design a GMsFEM for this nonlinear problem. In this way, we are able to construct an online-offline multiscale basis with respect to this nonlinearity.
%
%
We now outline the general procedure of the GMsFEM algorithm.

We begin briefly with some standard notation. The overall fine-scale model equations will be solved on a fine-grid using spaces $V_{h}, Q_{h}$ and $\hat{V}_{h}, \hat{Q}_{h}$, and will  be used for our reference solutions. 
We now introduce the coarse grid. 
%
%
Let $\mathcal{T}^H$ be a standard conforming partition of the computational domain $\Omega$ into finite elements. The fine-grid, ${\cal T}_{h}$ can be taken as a refinement of the coarse-grid. 
We refer to this partition as the coarse-grid and assume that each coarse element is partitioned into a connected union of fine grid blocks. 
We use $\{x_i\}_{i=1}^{N}$, where $N$ is the number of coarse nodes, to denote the vertices of the coarse mesh $\mathcal{T}^H$, and define the neighborhood of the node $x_i$ by
\[
\omega_i=\bigcup_{j}\left\{ K_j\in\mathcal{T}^H \, | \,  x_i\in \overline{K}_j\right\}.
\]
\begin{figure}[htb]
  \centering
  \includegraphics[width=0.6 \textwidth]{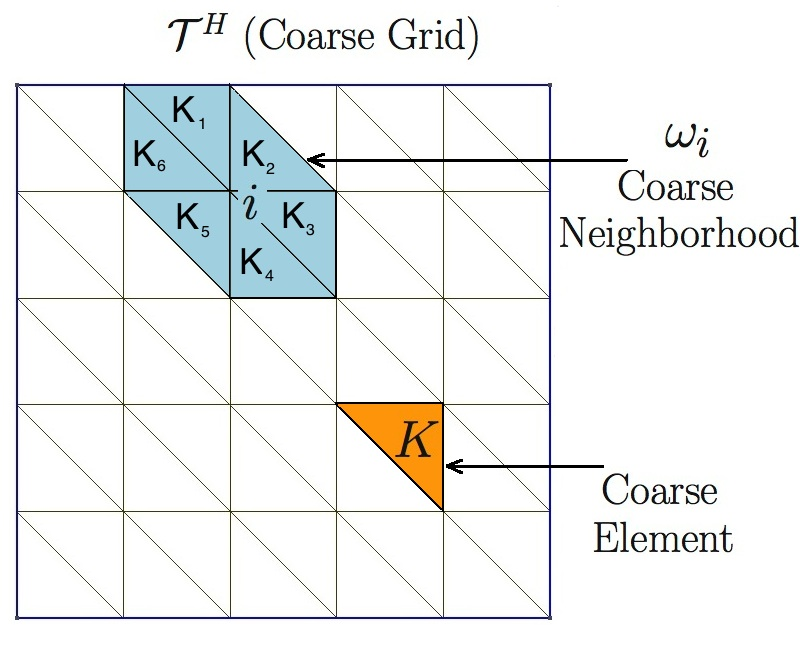}
  \caption{Illustration of a coarse neighborhood and coarse element}
  \label{schematic}
\end{figure}
See Figure~\ref{schematic} for an illustration of neighborhoods and elements subordinated to the coarse discretization. We emphasize that the use of $\omega_i$ is to denote a coarse neighborhood, and we use $K$ to denote a coarse element throughout the paper.

For global coupling we use the linearized continuous Galerkin (CG) formulation to find $(u^{n+1}_{j+1},p^{n+1}_{j+1}) \in (Q_{\text{on}},V_{\text{on}})$  such that 
\begin{subequations}
\label{CG_coupling}
\begin{eqnarray}
a^L( \mu_j ;u^{n+1}_{j+1},  v) + g(p^{n+1}_{j+1},  v) = 0,  \\
d \left( \frac{ u_{j+1}^{n+1} -  u^n}{\tau}, q \right)  
+ c \left( \frac{p^{n+1}_{j+1} - p^n}{\tau}, q \right)
+  b^L(\mu_j ; p^{n+1}_{j+1},  q)  = (f , q),  
\end{eqnarray}
\end{subequations}
where  $Q_{\text{on}}$ and $V_{\text{on}}$   denote the online spaces. The online spaces  are spanned by multiscale
 basis functions $\psi_m^{\omega_i, \text{on}}(x, \mu)$ and 
 $\varphi_{k}^{\omega_i, \text{on}}(x, \mu)$ for $n+1$ time step and $j$-th iteration, each of which is supported in $\omega_i$
\[
p(x, t)=\sum_{i,m} p_{m}^i(t) \psi_{k}^{\omega_i, \text{on}}(x, \mu), \quad
u(x, t)=\sum_{i,k} u_{k}^i(t) \varphi_{k}^{\omega_i, \text{on}}(x, \mu).
\]
The indexes $m, k$ represent the numbering of these multiscale  basis functions for pressure and displacements, respectively. 
Here the parameter $\mu$ represents the nonlinear dependence as in \eqref{CG_coupling}.  Recall that we may take $\mu_j=(u_{j},\nabla u_{j},p_{j},\nabla p_{j})$, from the previous time-step and will treat these variables as parametric values on each coarse patch. 
However, for simplicity we will suppose that the dependence is only on $(u,p)$ of this nonlinearity. 

 \vspace{.1in}

 {\bf Remark}:
Note, the derivative dependent problems, with nonlinear couplings of $(\nabla u ,\nabla p)$, may be handled.  
However, due to the oscillation in these quantities, these terms may not be well approximated by  constants on the coarse-grid level. Thus, we would need to have a more enriched parameter space than is utilized here. 

\vspace{.1in}

We now discuss further how we handle the parametrized nonlinearities. We assume that $u$ and $p$ are bounded above and below, i.e. $u \in [ u_{min}, u_{max} ] $ and $p \in [ p_{min}, p_{max}]$, where $(u_{min}, u_{max})$ and $(p_{min}, p_{max})$ are pre-defined constants.  These may be guessed initially based on initial data or a-priori estimates. 
The intervals $[u_{min}, u_{max} ] $ and $[ p_{min}, p_{max}]$ are divided into $N$ equal regions: 
$$u_{min} = u_0 < u_1 < ... < u_{N-1} < u_N= u_{max}, $$ and $$p_{min}= p_0 < p_1 < ... < p_{N-1} < p_N=p_{max}.$$
Clearly, if necessary these domains can be partitioned in different number of regions, but for simplicity we suppose they are equal in number. 
For the parameter $\mu_j$ we take average values of $u^{n+1}_j$ and $p^{n+1}_j$ in each coarse region $\omega_i$. For average of a function we will use the notation 
$$\bar{f}=\frac{1}{|\omega_i|}\int_{\omega_i} fdx.$$ 
More specifically, we use $\mu_j$ to represent the dependence of the solution on  $(\bar{u}^{n+1}_{j},\bar{p}^{n+1}_{j})$. 
The multiscale basis functions will be computed for a selected number of the parameter values $\mu_j$, $j=0,...,N$ at the offline stage and we will compute multiscale basis functions for each new value of $(\bar{u}^{n+1}_{j},\bar{p}^{n+1}_{j})$ for each $\omega_i$ at the online stage.

Boadly speaking, the GMsFEM algorithm consist of several steps:
\begin{itemize}
\item {\bf Offline computations:} 
\begin{enumerate}
\item Generate the coarse-grid, $\mathcal{T}^H$.
\item Construct the  snapshot space, used to compute an offline space, by solving many local problems on the fine-grid.
\item Construct  a small dimensional offline space by performing dimension reduction in the space of local snapshots.
\end{enumerate}
\item {\bf Online computations: } 
\begin{enumerate}
\item In each time step and nonlinear iteration for current value of $\mu_j$ in each $\omega_i$, we compute multiscale basis functions and construct online space by performing dimension reduction in the  offline space.
\item Use small dimensional online space to find the solution of a coarse-grid problem for any force term and/or boundary condition.
\end{enumerate}
\end{itemize}

We  construct multisclate basis functions for pressure and displacements separately. We begin by considering the pressure solve, then, the displacement solve.

\subsection{Multiscale basis functions for pressure}\label{pressuresolve}
In the offline computation, we first construct a snapshot space $Q_{\text{snap}}^{\omega}$. 
Construction of the snapshot space  involves solving the local problem for various choices of input parameters  and various boundary conditions. 
These local spatial fields are used then used construct the offline space 
and  the space consists of fields defined on a fine grid. There are a few options available when constructing the snapshot space and we will proceed with the two most natural ways.

\textbf{Snapshot Space 1:} First, we propose a  snapshot space generated by harmonic extensions of $b^L$. For simplicity, we will omit the index $i$ when there is no ambiguity. We thus define $\psi_{l,j}^{\omega, \text{snap}}$ such that 
\begin{equation} 
\label{harmonic_ex}
\begin{split}
b^L(\mu_j ; \psi_{l,j}^{\omega, \text{snap}} , q)&= 0 \quad \text{ in } \, \omega, \\
\psi_{l,j}^{\omega, \text{snap}}&=\delta_l^h(x) \quad \text{ on } \partial\omega.
\end{split}
\end{equation}
Here $\delta_l^h(x)$  are defined by
$\delta_l^h(x)=\delta_{l,k},\,\forall l \in \textsl{J}_{h}(\omega)$, where $\textsl{J}_{h}(\omega)$ denotes the fine-grid boundary node on $\partial\omega$.  This is done for each fixed parameter $\mu_j,$ $j = 0,...,N.$

\textbf{Snapshot Space 2:} Alternatively, we may use local fine-scale space basis functions within a coarse region and construct local snapshots by solving the following eigenvalue problem with natural boundary conditions
\begin{equation} 
\label{snapshot2}
\begin{split}
B(\mu_j) \psi_{l,j}^{\omega, \text{snap}} &= \lambda_{l,j}^{\omega, \text{snap}} M(\mu_j) \psi_k^{\omega, \text{snap}}, \quad \text{ in } \, \omega.
\end{split}
\end{equation}
Where  
\[
B_{ij}(\mu_j) = \int_{\Omega} \left( K(x, \mu_j) \nabla \phi_i , \nabla \phi_j \right)\, dx, 
\quad
M_{ij}(\mu_j) = \int_{\Omega}  K(x, \mu_j)  \phi_i \phi_j \, dx,
\]
   $\phi_i$ are the standard fine-scale basis functions,
and 
for each fixed parameter values $\mu_j$, $j = 0,...,N$. 

Let $l_i$ be the number of functions in the snapshot space in the region $\omega$, and define
\[
Q^{\omega}_{\text{snap}} = \text{span}\{ \psi_{l,j}^{ \text{snap}}:\quad 1\leq l \leq l_i, \quad 0 \leq j \leq N \},
\]
for each coarse subdomain $\omega$.
We reorder the snapshot functions using a single index to create the matrix
$$
R^p_{\text{snap}} = \left[ \psi_{1}^{\text{snap}}, \ldots, \psi_{M_{\text{snap}}}^{\text{snap}} \right],
$$
where $M_{\text{snap}}$ denotes the total number of functions to keep in the snapshot  construction. 

To construct the offline space $Q_{\text{off}}$, we perform a dimension reduction of the space of snapshots by using an auxiliary spectral decomposition. More precisely,  we solve the  eigenvalue problem in the space of snapshots:
\begin{equation} 
\label{offeig}
B^{\text{off}} \Psi_k^{\text{off}} = \lambda_k^{\text{off}} M^{\text{off}} \Psi_k^{\text{off}},
\end{equation}
where
\[
\begin{split}
B^{\text{off}} &= \int_{\omega} \left( \overline{K}(x) \nabla \phi^{\text{snap}}_i , \nabla \phi^{\text{snap}}_j \right)\, dx 
= (R^p_{\text{snap}})^T \overline{B} R^p_{\text{snap}},
\\
M^{\text{off}} &= \int_{\omega}  \overline{K}(x)  \phi^{\text{snap}}_i \phi^{\text{snap}}_j \, dx 
= (R^p_{\text{snap}})^T \overline{M} R^p_{\text{snap}}.
\end{split}
\]
Here $$\overline{K}(x) = \sum_{j=1}^{N} t_j K(x, \mu_j),$$ is independent of $\mu_j$ and $t_j$ are prescribed non-negative weights.  
The main objective is to use the offline space to accurately construct a set of multiscale basis functions for each $\mu_j$ in the online stage. At the offline stage the bilinear forms are chosen to be parameter-independent, such that there is no need to reconstruct the offline space for each $\mu_j$. 

We then  choose the smallest $N^{\omega,p}_{\text{off}}$ eigenvalues from Eq.~\eqref{offeig} and form the corresponding eigenvectors in the space of snapshots by setting
$ \psi_k^{\text{off}} = \sum_{j=1}^{M_{\text{snap}}} \Psi_{kj}^{\text{off}} \psi_j^{\text{snap}}$,  for $k=1,\ldots, N^{\omega,p}_{\text{off}}$, where $\Psi_{kj}^{\text{off}}$ are the coordinates of the vector $\psi_{k}^{\text{off}}$. We denote the span of this reduced space as $Q_{\text{off}}^{\omega}$.

At the online stage, for a given parameter value $\mu$, multiscale basis functions are computed based on each local coarse region $\omega_i$.  The associated online space $Q_{\text{on}}^{\omega}(\mu)$  is the small dimensional subspace of the offline space.
 In particular, we seek  a subspace of the offline space that can approximate any element of the offline space in an appropriate sense.  
 In the the online stage the bilinear forms are chosen to be parameter-dependent and we use following eigenvalue problem
 \begin{equation} 
\label{oneig}
B^{\text{on}} \Psi_k^{\text{on}} = \lambda_k^{\text{on}} M^{\text{on}} \Psi_k^{\text{on}},
\end{equation}
 where 
 \[
\begin{split}
B^{\text{on}} &= \int_{\omega} \left( K(x, \mu) \nabla \phi^{\text{off}}_i , \nabla \phi^{\text{off}}_j \right)\, dx 
= (R^p_{\text{off}})^T B R^p_{\text{off}},
\\
M^{\text{on}} &= \int_{\omega} K(x, \mu)  \phi^{\text{off}}_i \phi^{\text{off}}_j \, dx 
= (R^p_{\text{off}})^T M R^p_{\text{off}}.
\end{split}
\]
Here $B$ and $M$ are the fine scale matrices corresponding to the stiffness and mass matrices for given $\mu$ and 
$$
R^p_{\text{off}} = \left[ \psi_{1}^{\text{off}}, \ldots, \psi_{N^{\omega,p}_{\text{off}}}^{\text{off}} \right].
$$

Finally, we  multiply the partition of unity functions $ \chi_i $ by the eigenfunctions in the online space $Q_{\text{on}}^{\omega_i}$ to construct the resulting basis functions
\begin{equation} 
\label{cgbasis}
\psi_{i,k} = \chi_i \psi_k^{\omega_i, \text{on}} \quad \text{for} \, \, \,
1 \leq i \leq N_c \, \, \,  \text{and} \, \, \, 1 \leq k \leq N_{\text{on}}^{\omega_i,p},
\end{equation}
where $\psi_k^{\text{on}} = \sum_{j=1}^{l_i} \Psi_{kj}^{\text{on}} \psi_k^{\text{off}}$, $\chi_i$ is the standard linear partition of unity functions and the $N_{\text{on}}^{\omega_i,p}$ denotes the number of online eigenvectors that are chosen for each coarse node $i$. We note that the construction in Eq.~\eqref{cgbasis} yields  continuous basis functions due to the multiplication of offline eigenvectors with the initial (continuous) partition of unity. Next, we define the online space as
\begin{equation} \label{cgspace}
Q_{\text{on}}  = \text{span} \{ \psi_{i,k} : \,  \, 1 \leq i \leq N_c \, \, \,  \text{and} \, \, \, 1 \leq k \leq N_{\text{on}}^{\omega_i,p}  \}.
\end{equation}
Using a single index notation, we may write $Q_{\text{on}} = \text{span} \{ \psi_{i} \}_{i=1}^{N^p_c}$, where $N^p_c =\sum_{i=1}^{N_c}N_{\text{on}}^{\omega_{i},p}$  denotes the total number of basis functions in the spaces $Q^{\omega_{i}}_{\text{on}}$ and $N_c$ is number of coarse mesh nodes. 

Denote the matrix  
\[
R_p = \left[ \psi_1 , \ldots, \psi_{N^p_c} \right]^T,
\] 
where $\psi_i$ are used to denote the nodal values of each basis function defined on the fine grid.

\subsection{Multiscale basis functions for displacements}

For construction of multiscale basis functions for displacements we use similar algorithm that we used for pressure. 
We first construct a snapshot space $V_{\text{snap}}^{\omega}$ for each parameter $\mu_j$. Again, as with pressure we give two possible snapshot space choices. 

\textbf{Snapshot Space 1:} As our first possible snapshot space we propose the harmonic extension using $a^L$. We define $\varphi_{l,j}^{\omega, \text{snap}}$ as the solution to 
\begin{equation} 
\label{harmonic_ex2}
\begin{split}
a^L( \mu_j ; \varphi_{l,j}^{\omega, \text{snap}} , v) &= 0 \quad \text{in } \, \omega, \\
\varphi_{l, j}^{\omega, \text{snap}}&=\delta_l^h(x), \quad \text{ on }\partial\omega.
\end{split}
\end{equation}
Again,  
$\delta_l^h(x)=\delta_{l,k},\,\forall l \in \textsl{J}_{h}(\omega)$, and 
for each fixed parameter values $\mu_j$, $j = 0,...,N$. 

\textbf{Snapshot Space 2:} We could also use the method based on solving an eigenvalue problem with natural boundary conditions given by 
\begin{equation} 
\label{usnapshot2}
\begin{split}
A(\mu_j) \Phi_{l,j}^{\omega, \text{snap}} &= \lambda_{l,j}^{\omega, \text{snap}} N(\mu_j) \Phi_k^{\omega, \text{snap}}, \quad \text{ in } \, \omega,
\end{split}
\end{equation}
Where 
\[
A_{ij}(\mu_j) = \int_{\Omega} (C (x; \mu_j)\varepsilon(\varphi_i) , \varepsilon(\varphi_j) ) \, dx, 
\quad
N_{ij}(\mu_j) = \int_{\Omega}  m (x; \mu_j) \varphi_i \varphi_j \, dx,
\]
and, in the case of linear elasticity   $ m (x; \mu_j) = (\lambda + 2 \mu_{e})$. In a more complicated relation $ m (x; \mu_j) $ is related to the lower order operators \cite{ElasticGMsFEM}.  Again,  $\varphi_i$ are the standard fine-scale basis functions, and his is done for each fixed parameter values $\mu_j$, $j = 0,...,N$.

Define
\[
V^{\omega}_{\text{snap}} = \text{span}\{ \Phi_{l,j}^{ \text{snap}}:\quad 1\leq l \leq l_i, \quad 0 \leq j \leq N \},
\]
for each coarse subdomain $\omega$. 
We denote the corresponding matrix of snapshot functions, again with similar notation, to be 
$$
R^u_{\text{snap}} = \left[ \Phi_{1}^{\text{snap}}, \ldots, \Phi_{N_{\text{snap}}}^{\text{snap}} \right].
$$
where $N_{\text{snap}}$ denotes the total number of functions to keep in the snapshot  construction.

Again, we perform a dimension reduction of the space of snapshots by using an auxiliary spectral decomposition. We solve the parameter-independent eigenvalue problem in the space of snapshots
\begin{equation} 
\label{offeig2}
A^{\text{off}} \Phi_k^{\text{off}} = \lambda_k^{\text{off}} N^{\text{off}} \Phi_k^{\text{off}},
\end{equation}
where
\[
A^{\text{off}} = (R^u_{\text{snap}})^T A R^u_{\text{snap}},
\quad
N^{\text{off}} = ( R^u_{\text{snap}})^T N R^u_{\text{snap}},
\]
where $A$ and $N$ denote fine scale matrices 
\[
A_{mn}
= \int_{\omega} ( \overline{C} (x)  {\varepsilon}( \varphi_m) , {\varepsilon}( \varphi_n ) \, dx,
\quad
 N_{mn}
= \int_{\omega} \overline{m} (x) \varphi_m \cdot  \varphi_n  \, dx .
\]
Here, $\varphi_i$ are  fine-scale basis functions. Further, we have 
$$\overline{C}(x) = \sum_{j=1}^{N} t_j C(x, \mu_j), \quad \overline{m}(x) = \sum_{j=1}^{N} t_j m(x;\mu_j)$$
 is independent of $\mu_j$ and $t_j$ are prescribed non-negative weights.  
Recall, the main objective is to use the offline space to accurately construct a set of multiscale basis functions for each $\mu_j$ in the online stage. 
As before for the fluids flow module, at the offline stage of the mechanics the bilinear forms are chosen to be parameter-independent, such that there is no need to reconstruct the offline space for each $\mu_j$.

We then  choose the smallest $N^{\omega,u}_{\text{off}}$ eigenvalues from Eq.~\eqref{offeig2} and form the corresponding eigenvectors in the space of snapshots by setting
$\varphi_k^{\text{off}} = \sum_{j=1}^{l_i} \Phi_{kj}^{\text{off}} \Phi_j^{\text{snap}}$,  for $k=1,\ldots, N^{\omega,u}_{\text{off}}$, where $\Phi_{kj}^{\text{off}}$ are the coordinates of the vector $\varphi_{k}^{\text{off}}$. We denote the span of this reduced space as $V_{\text{off}}^{\omega}$ and denote
$$
R^u_{\text{off}} = \left[ \varphi_{1}^{\text{off}}, \ldots, \varphi_{N^{\omega,u}_{\text{off}}}^{\text{off}} \right].
$$

At the online stage, we use following parameter-dependent eigenvalue problem
 \begin{equation} 
\label{oneig2}
A^{\text{on}} (\mu) \Phi_k^{\text{on}} = \lambda_k^{\text{on}} N^{\text{on}} (\mu) \Phi_k^{\text{on}},
\end{equation}
 where 
 \[
\begin{split}
A^{\text{on}} (\mu) &= \int_{\omega} \left( C (x; \mu)  {\varepsilon}( \varphi_m^{\text{off}}) , {\varepsilon}( \varphi_n^{\text{off}})\right)\, dx 
= (R^u_{\text{off}})^T A R^u_{\text{off}},
\\
N^{\text{on}}(\mu)  &= \int_{\omega} m(x, \mu)  \varphi^{\text{off}}_m \varphi^{\text{off}}_n \, dx 
= (R^u_{\text{off}})^T N R^u_{\text{off}}.
\end{split}
\]

Finally, we  multiply the linear partition of unity functions $\xi_i$ by the eigenfunctions in the online space $V_{\text{on}}^{\omega_i}$ to construct the resulting basis functions
\begin{equation} 
\label{cgbasis.mechanics}
\varphi_{i,k} = \xi_i \varphi_k^{\omega_i, \text{on}} \quad \text{for} \, \, \,
1 \leq i \leq N_c \, \, \,  \text{and} \, \, \, 1 \leq k \leq N_{\text{on}}^{\omega_i,u},
\end{equation}
where $\varphi_k^{\text{on}} = \sum_{j=1}^{l_i} \Phi_{kj}^{\text{on}} \varphi_j^{\text{off}}$ and $N_{\text{on}}^{\omega_i,u}$ denotes the number of online eigenvectors that are chosen for each coarse node $i$. 
Next, we define the online space as
\begin{equation} 
\label{cgspace.mechanics}
V_{\text{on}}  = \text{span} \{ \varphi_{i,k} : \,  \, 1 \leq i \leq N_c \, \, \,  \text{and} \, \, \, 1 \leq k \leq N_{\text{on}}^{\omega_i}  \}.
\end{equation}
Using a single index notation, we may write $V_{\text{on}} = \text{span} \{ \varphi_{i} \}_{i=1}^{N^u_c}$, where $N^u_c =\sum_{i=1}^{N_c}N_{\text{on}}^{\omega_{i},u}$
denotes the total number of basis functions in the space $V^{\omega_{i}}_{\text{on}}$.

And after construction $V_{\text{on}}$ we 
denote the matrix  
\[
R_u = \left[ \varphi_1 , \ldots, \varphi_{N^u_c} \right]^T,
\] 
where $\varphi_i$ are used to denote the nodal values of each basis function defined on the fine grid.

\subsection{Global coupling}
Now that we have constructed the online spaces for both the fluid and mechanics we now can use this parametrized basis at the global level. Indeed, 
for global coupling we use system of equations \eqref{CG_coupling}  to find $(u^{n+1}_{j+1},p^{n+1}_{j+1}) \in (Q_{\text{on}},V_{\text{on}})$, where
\[
Q_{\text{on}} = \text{span} \{ \psi_{i} \}_{i=1}^{N^p_c},
\quad \text{and} \quad 
V_{\text{on}} = \text{span} \{ \varphi_{i} \}_{i=1}^{N^u_c}.
\]
Using the matrices
\[
R_p = \left[ \psi_1 , \ldots, \psi_{N^p_c} \right]^T,
\quad \text{and} \quad 
R_u = \left[ \varphi_1 , \ldots, \varphi_{N^u_c} \right]^T,
\] 
we may write matrix analogue for the  variational for \eqref{CG_coupling} that  will be used for calculation of multiscale solution $(u_{j+1}^{ms, n+1}, p_{j+1}^{ms, n+1}).$
Writing \eqref{CG_coupling} in matrix form,  using the notation in \eqref{eq:mmd}, in the online basis we have 
\begin{eqnarray}
\label{eq:CG-matrixform}
R_u A(\mu_j) R^T_u u_{j+1}^{ms, n+1}  + \alpha R_u G R^T_p p_{j+1}^{ms, n+1} = 0,  \\
R_p (S + \tau B(\mu_j)  )R^T_p  p_{j+1}^{ms, n+1} + \alpha R_p  D  R^T_u u_{j+1}^{ms, n+1}  = \tau R_p F + R_p S R^T_p p^n.
\end{eqnarray}
We also note that matrices $R_p$ and $R_u$ may be analogously used in order to project coarse-scale solutions onto the fine-grid 
$$p^{n+1}_{j+1} =  R_p^T p_{j+1}^{ms, n+1}, \quad  u^{n+1}_{j+1} =  R_u^T u_{j+1}^{ms, n+1}.$$

\section{Numerical Examples}\label{Numerics}
In this section, we present numerical examples to  demonstrate the  performance of the GMsFEM  for computing the solution of the nonlinear poroelasticity problem in heterogenous domains and complex nonlinear dependence on permeability and elastic properties. 
We use fully coupled scheme for approximation by time with Picard iteration to linearize the nonlinearity.
We will implement a single complicated geometry with contrasting parameter values. Indeed, as noted before, there are many possible nonlinear relations, but here we take a an exponential pressure relationship with the permeability.
We present the errors with varying number of multiscale basis functions and over time for linear and nonlinear case with  parameters. 

\begin{figure}[htb]
\begin{center}
\includegraphics[width=0.35\linewidth]{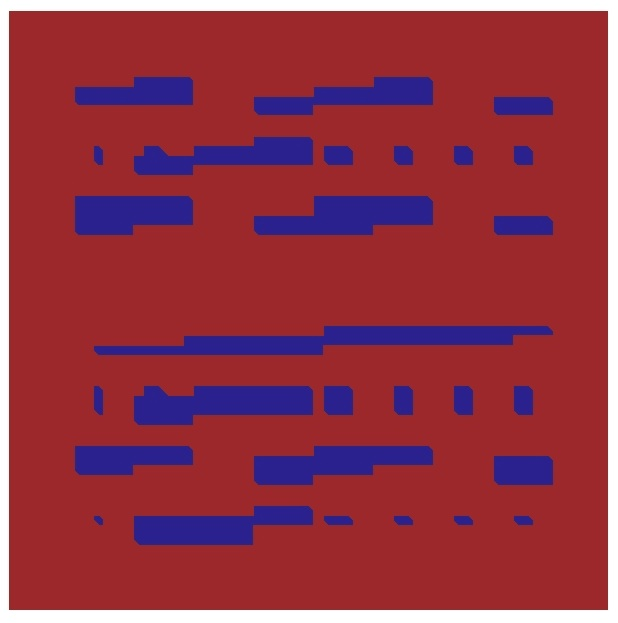} 
\caption{Coefficients subdomains. Red is the subdomain 1 and blue is the subdomain 2}
\label{fig:koeff}
\end{center}
\end{figure}

We proceed as in \cite{BrownVasilyeva}, and 
we take the computational domain $\Omega$ as a unit square $[0,1]^2$, and set the source term $f = 0$ in \eqref{eq:main}.
We utilize  heterogeneous coefficients that have different values in  two subdomains.  We denote each region as subdomain 1 and 2, $\Omega_{1},\Omega_{2}$, respectively.
We use following coefficients: for the Biot modulus we take $M_1 = 1.0, M_2 = 10$ in each respective numbered subdomain $\Omega_{i}$. For permeability we take a linear $K$ and nonlinear relation $K(p)$. 
More specifically, for the linear regime we have 
%
\begin{align}\label{linearK}
K =
\begin{cases}
\exp(1) & \text{ in  } \Omega_{1}, \\
\exp(10) & \text{ in } \Omega_{2}.
\end{cases}
\end{align}
For nonlinear case  we use a permeability that depends on pressure $p$
\begin{align}\label{nonlinearK}
K(p) =
\begin{cases}
\exp(p) & \text{ in  } \Omega_{1}, \\
\exp(10p) & \text{ in } \Omega_{2}.
\end{cases}
\end{align}
For fluid-solid coupling constant we have $\alpha = 0.9$. 
For the elastic properties we use following coefficients: elastic modulus is given by  $E_1 = 10, E_2 = 1$ in each respective subdomain $\Omega_{i}$, the Poisson's ratio is $\eta = 0.22$, and these can be related to the parameters $\mu_{i}$ and $\lambda_{i}$, for $i=1,2,$ via the relation
\[
\mu_{i} = \frac{E_{i}}{2 (1 + \eta)}, \quad
\lambda_{i} = \frac{E_{i} \eta}{(1+ \eta) ( 1- 2 \eta)},
\]
in each subdomain. The subdomains for coefficients shown in Figure \ref{fig:koeff},  where the background media in red is the subdomain 1, $\Omega_{1},$ and  isolated particles  and strips in blue are the subdomain 2, $\Omega_{2}.$

As we have chosen $f=0$ we must use boundary conditions to force flow and mechanics. In these tests, as in \cite{BrownVasilyeva}, we use following boundary conditions: 
\[
p = p_1, \quad x \in \Gamma_T, \quad
p = p_0, \quad x \in \Gamma_B, \quad
\frac{\partial p}{\partial n} = 0, \quad x \in \Gamma_L \cup \Gamma_R,
\] 
and
\[
u_{x} = 0, \quad \frac{\partial u_{y}}{\partial y} = 0, \quad x \in \Gamma_L, \quad
\frac{\partial u_x}{\partial x} = 0, \quad u_{y} = 0, \quad x \in \Gamma_B, 
\]
and finally,
\[
\frac{\partial u_x}{\partial x} = 0, \quad \frac{\partial u_y}{\partial y} = 0, \quad x \in \Gamma_T \cup \Gamma_R.
\]
Here $\Gamma_L$ and $\Gamma_R$ are left and right boundaries, $\Gamma_T$ and $\Gamma_B$ are top and bottom boundaries respectively. We set $p_0 = 0$ and $p_1 = 1$ to drive the flow, and thus, the mechanics. 

\begin{figure}[htb]
\begin{center}
\includegraphics[width=1\linewidth]{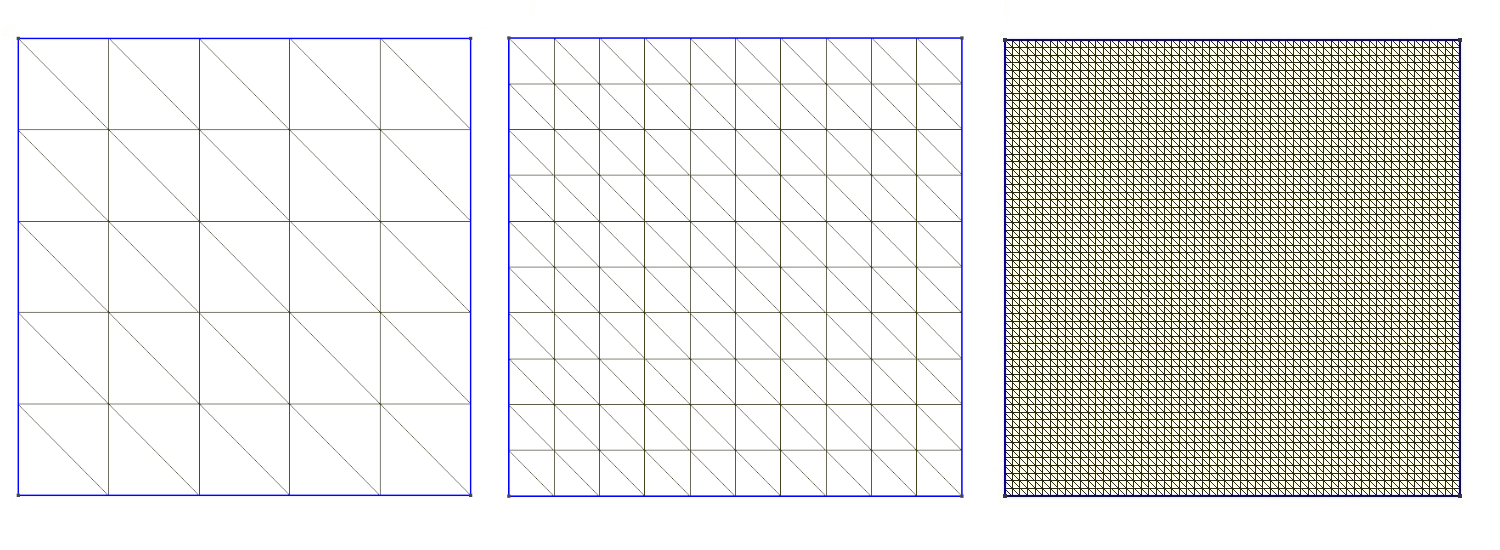} 
\caption{Two coarse grids and fine grid. Left: coarse grid with 36 nodes. Middle: coarse grid with 121 nodes. Right: fine grid with 3721 nodes.}
\label{fig:domain}
\end{center}
\end{figure}

In Figure \ref{fig:domain} we show the two coarse grids and fine grid. The first coarse grid  consists of 36 nodes and 50 triangle cells, the second coarse grid contains 121 nodes and 200 triangle cells, and the fine mesh consists of 3721 nodes and 7200 triangle cells. 
The number of time steps  is $M_{T}=10$ and the maximal time being set at $T_{max} = 0.055$. 
As an  initial condition for pressure we use $p = p_0=0$.
For the nonlinear solve  we use Picard iteration for linearization and terminate the iterative loop when $|| p_f - p_{ms}||_{L_2(\Omega)} \leq \delta$, $\delta = 10^{-5}$.

\begin{figure}[htb]
\begin{center}
\includegraphics[width=0.8\linewidth]{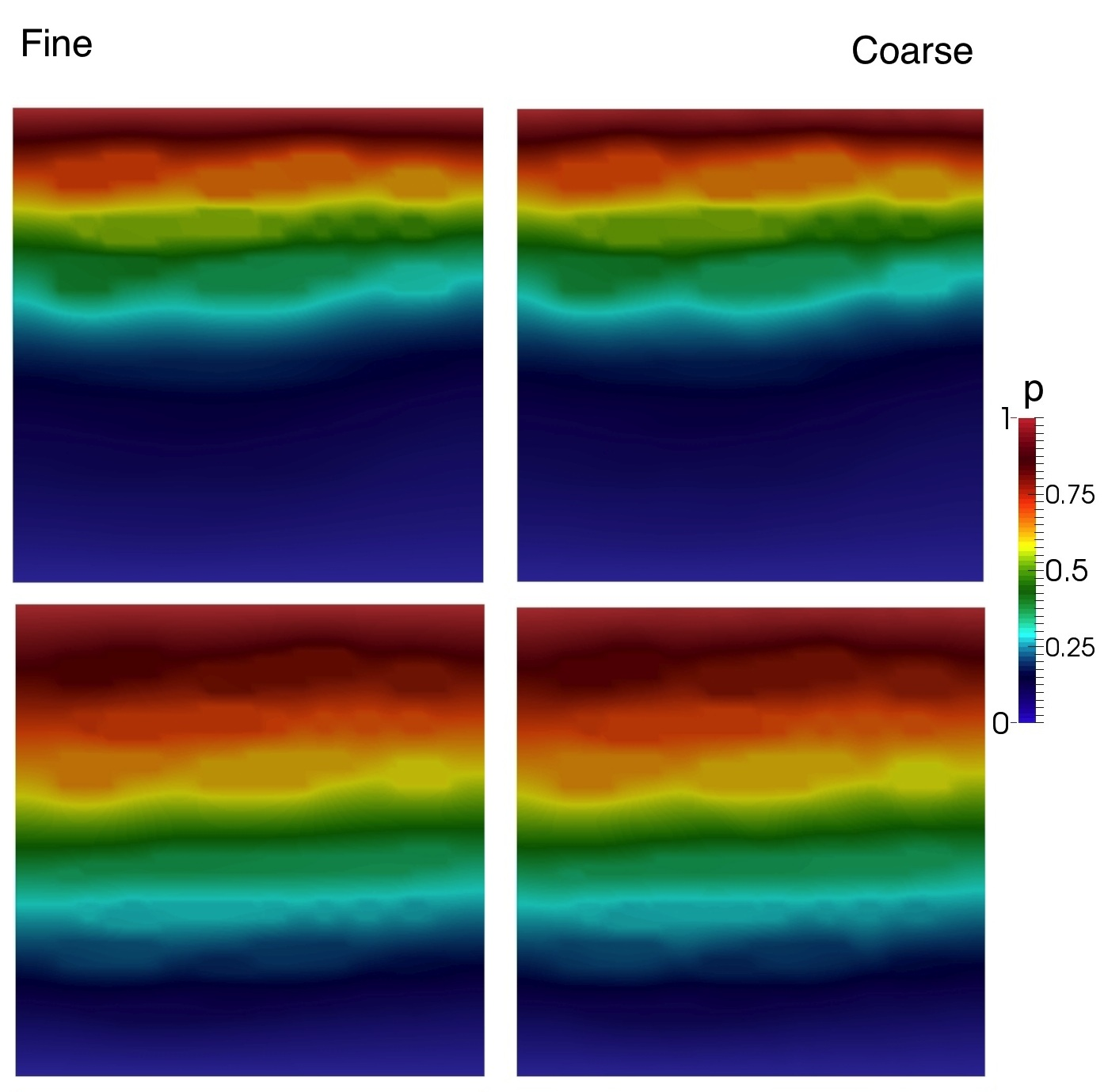} 
\caption{The fine-scale and coarse-scale solutions of the pressure distribution for $T = 0.02$ and $0.055$ (from top to bottom) for nonlinear case. The dimension of the fine-scale solution is 11163 and the dimension of the coarse space is 864.}
\label{fig:p}
\end{center}
\end{figure}

\begin{figure}[htb]
\begin{center}
\includegraphics[width=0.8\linewidth]{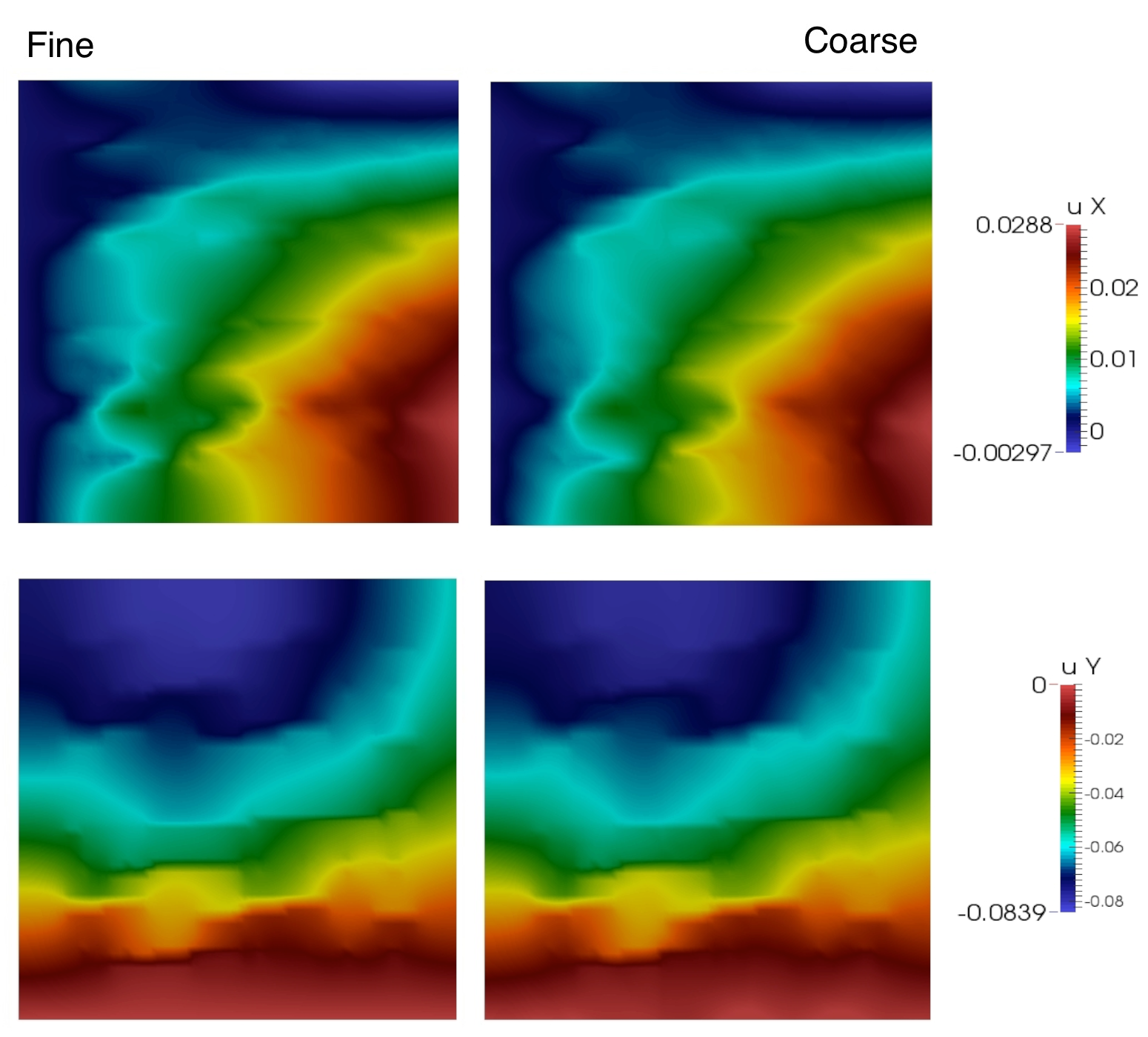} 
\caption{The fine-scale and coarse-scale solutions of the displacements $u_x$ and $u_y$ for nonlinear case. The dimension of the fine-scale solution is 11163 and the dimension of the coarse space is 864.}
\label{fig:u}
\end{center}
\end{figure}

The reference solution computed by using a  standard FEM (linear basis functions for pressure and displacements) on the fine grid, Picard type linearization, and   using a fully coupled time-splitting scheme.
The pressure and the displacement fields on the fine-grid are presented on the left column of  Figure \ref{fig:p}  and Figure \ref{fig:u}.

The errors will be measured in relative weighted $L^2$ and relative weighted $H^1$ norm for pressure
\begin{eqnarray*}
\norm{\varepsilon_p}_{L^2(\Omega)} &= 
\frac{   \left(  \int_{\Omega} K(p_f) (p_f - p_{ms})^2 dx \right)^{1/2}}{  \left(  \int_{\Omega}K(p_f) p_f^2 dx \right)^{1/2} } , \\
\seminorm{\varepsilon_p}_{H^1(\Omega)} &= \frac{ \left(  \int_{\Omega} \left( K(p_f) \nabla (p_f - p_{ms}), \nabla  (p_f - p_{ms}) \right) dx \right)^{1/2}  }{ \left(  \int_{\Omega} \left( K(p_f)\nabla p_f, \nabla p_f \right) dx \right)^{1/2}},
\end{eqnarray*}
and for displacements, due to the linearity in our Elasticity in this example, we have 
\begin{eqnarray*}
\norm{\varepsilon_u}_{L^2(\Omega)} &= \frac{  \left(  \int_{\Omega} (\lambda + 2 \mu) (u_f - u_{ms}, u_f - u_{ms}) dx \right)^{1/2} }{ \left(  \int_{\Omega} (\lambda + 2 \mu) (u_f , u_f) dx \right)^{1/2} }, \\
\seminorm{\varepsilon_u}_{H^1(\Omega)} &= \frac{   \left(  \int_{\Omega} \left( \sigma (u_f - u_{ms}), \varepsilon(u_f - u_{ms}) \right) dx  \right)^{1/2} }{ \left(  \int_{\Omega} \left( \sigma (u_f), \varepsilon(u_f) \right) dx  \right)^{1/2}  }.
\end{eqnarray*}
Here $(u_f,p_f)$ and $(u_{ms},p_{ms})$ are  fine-scale and coarse-scale using GMsFEM solutions, respectively for pressure and displacements.

In our examples, the nonlinearity resides in the pressure solves. Therefore, we will use the nonlinear parameter dependence approach in Section \ref{pressuresolve}. For our Elasticity basis construction, we may remain in the linear algorithmic approach to construct the online basis. 
In general, for simulation using GMsFEM we first generate a snapshot space using first choice (\eqref{harmonic_ex}, snapshot space 1) or second choice (\eqref{snapshot2}, snapshot space 2), then we use a spectral decomposition to obtain the offline space, and similarly  to obtain the online space. 
For each time step and nonlinear Picard iteration we update the online space for pressure and solve the equation \eqref{oneig} utilizing the  previously computed solution $p_{j}^{n+1}$.
For construction the snapshot space 2 we choose a specified number of eigenfunctions $l_i = 16$ for all $\omega_i$.
We select the range of solutions $p_{min} = 0$ and $p_{max} = 1$ and divide the domain $[p_{min}, p_{max}]$ into $N$ equally spaced subdomains to obtain $N+1$ discrete points $p_0, ..., p_N$. For simulation we use $N = 20$. 

Recall, we will use a few multiscale basis functions per each coarse node $\omega_{i}$, and these number of coarse basis defines the problem size (dimension of online spaces, $Q_{\text{on}}$ and $V_{\text{on}}$).  
We suppose that in each patch $\omega_{i}$ we take the same number of multiscale basis functions for pressure, $N^p_{\text{on}}=N^{\omega_{i},p}_{\text{on}}$, for all $\omega_{i}$.
Similarly for displacements we take  $N^u_{\text{on}}=N^{\omega_{i},u}_{\text{on}}$, for all $\omega_{i}$.
Varying the basis functions in both pressure and displacement multiscale spaces we record the errors at the final time. 
We note that the size of online space and the associated solution accuracy will depend on the number of eigenvectors ($N^p_{\text{on}}$ and $N^u_{\text{on}}$) that we keep in the online space construction.

We begin first with the purely linear case with $K$ given by \eqref{linearK}.
In Tables \ref{tab:c1-5-lin} and  \ref{tab:c1-10-lin}, we present the  relative weighted $L^2$ and $H^1$ errors for linear case of the coefficients in geometry Figure \ref{fig:koeff}
using the fully coupled time scheme on two coarse grids. In Table  \ref{tab:c1-5-lin} we have a coarse-grid of 36 nodes and in Table \ref{tab:c1-10-lin} we have a refined coarse grid with 121 nodes. We compare these to a fine-scale solution space with dimension 11163.
In these tables, $N_{\text{on}}^{p}$ and $N_{\text{on}}^{u}$ are number of multiscale basis functions for each neighborhoods, the second column show the dimension of the online space, the next two columns present the relative weighted $L^2$ and $H^1$ errors for pressure and last two columns show the relative weighted $L^2$ and $H^1$ errors for displacements.
We note that as the dimension of the online space increases, because we keep more eigenfunctions $N_{\text{on}}^{p}$, $N_{\text{on}}^{u}$
 in the space construction. We note for the less refined coarse-grid with 36 nodes the relative weighted $L^2$ errors decrease from 36.5\% to 0.07\% for pressure and from 24.3\% to 0.5\% for displacements and relative weighted $H^1$ errors decrease from 99.0\% to 2.7\% for pressure and from 37.7\% to 3.4\% for displacements. 
 We note for the  refined coarse-grid with 121 nodes the relative weighted $L^2$ errors decrease from 14.1\% to 0.01\% for pressure and from 26.9\% to 0.1\% for displacements and relative weighted $H^1$ errors decrease from 82.0\% to 1.6\% for pressure and from 36.1\% to 2.5\% for displacements. 
 We note that in this example, refining the coarse-grid is not as advantageous to more local basis functions per grid-block. 
Indeed,  with the less refined coarse-grid of 36 nodes and $N_{\text{on}}^{p}=N_{\text{on}}^{u}=12$ gives a very good percentage error for a space of dimension 1296 when compared to the more refined coarse-grid of 121 nodes and less eigenvectors $N_{\text{on}}^{p}=N_{\text{on}}^{u}=4$ with space of dimension  1452. 
 %
 %

\begin{table}[htp]
\begin{center}
\begin{tabular}[hp]{|c|c|cc|cc|}
\hline
 & & \multicolumn{2}{|c|}{Pressure errors, $\varepsilon_p$}
 & \multicolumn{2}{|c|}{Displacements errors, $\varepsilon_u$} \\
$N_{\text{on}}^p$  & dim($Q_{\text{on}},V_{\text{on}}$) & $L^2$ & $H^1$ & $L^2$ & $H^1$\\
\hline \hline
\multicolumn{6}{|c|}{$N_{\text{on}}^u = 4$}  \\
\hline
2     & 360 & 0.365	& 0.990   & 0.243  &  0.377  \\
4     & 432 & 0.057	& 0.435   & 0.238  &  0.370  \\
\hline 
\multicolumn{6}{|c|}{$N_{\text{on}}^u = 8$}  \\
\hline
2     & 648 & 0.365	& 0.990	& 0.108	& 0.207 \\
4     & 720 & 0.057	& 0.435	& 0.045	& 0.077 \\
8     & 864 & 0.001	& 0.059	& 0.017	& 0.072 \\
\hline 
\multicolumn{6}{|c|}{$N_{\text{on}}^u= 12$}  \\
\hline
2     & 936 & 0.365		& 0.990	& 0.111   &  0.199 \\
4     & 1008 & 0.057	& 0.435	& 0.042   &  0.045 \\
8     & 1152 & 0.001	& 0.059	& 0.007   &  0.034 \\
12   & 1296 & 0.0007	& 0.027	& 0.005   &  0.034 \\
\hline 
\end{tabular}
\end{center}
\caption{Numerical results for linear problem for coarse mesh with 36 nodes.}
\label{tab:c1-5-lin}
\end{table}

\begin{table}[htp]
\begin{center}
\begin{tabular}[hp]{|c|c|cc|cc|}
\hline
 & & \multicolumn{2}{|c|}{Pressure errors, $\varepsilon_p$}
 & \multicolumn{2}{|c|}{Displacements errors, $\varepsilon_u$} \\
$N_{\text{on}}^p$  & dim($Q_{\text{on}},V_{\text{on}}$) & $L^2$ & $H^1$ & $L^2$ & $H^1$\\
\hline \hline
\multicolumn{6}{|c|}{$N_{\text{on}}^u = 4$}  \\
\hline
2     & 1210 & 0.141 &  0.827 & 0.269 &  0.361 \\
4     & 1452 & 0.007 &  0.132 & 0.240 &  0.352 \\
\hline 
\multicolumn{6}{|c|}{$N_{\text{on}}^u = 8$}  \\
\hline
2     & 2178 & 0.141 &  0.827 & 0.069 &  0.095 \\
4     & 2420 & 0.007 &  0.132 & 0.024 &  0.063 \\
8     & 1904 & 0.001 &  0.042 & 0.015 &  0.062 \\
\hline 
\multicolumn{6}{|c|}{$N_{\text{on}}^u= 12$}  \\
\hline
2     & 3148 & 0.141 	&  0.827 & 0.059 & 0.076 \\
4     & 3388 & 0.007 	&  0.132 & 0.011 & 0.027 \\
8     & 3872 & 0.001 	&  0.042 & 0.003 & 0.025 \\
12   & 4356 & 0.0001 &  0.016 & 0.001 & 0.025 \\
\hline 
\end{tabular}
\end{center}
\caption{Numerical results for linear problem for coarse mesh with 121 nodes.}
\label{tab:c1-10-lin}
\end{table}

In a similar setting, we consider the nonlinear case of the coefficient with $K(p)$ given by \eqref{nonlinearK}. Here we will explore the different snapshot spaces available for us in the nonlinear algorithm. Again as in the linear case we use two coarse-grids and implement this with a fully coupled time scheme and use Picard iterations for the nonlinearity. 
We present the results in Table \ref{tab:c1-1} for snapshot space 1, the errors are very similar in magnitude when compared to the corresponding linear case. 
In the left side of Table \ref{tab:c1-1} we present the errors for 36 nodes in the coarse-grid.
The relative weighted $L^2$ errors decrease from 8.1\% to 0.09\% for pressure and from 30.4\% to 0.5\% for displacements  and relative weighted $H^1$ errors decrease from 60.9\% to 4.7\% for pressure and from 38.0\% to 3.4\% for displacements. 
In the right side of Table \ref{tab:c1-1} we present the errors for 121 nodes in the coarse-grid.
The relative weighted $L^2$ errors decrease from 4.8\% to 0.02\% for pressure and from 26.4\% to 0.1\% for displacements  and relative weighted $H^1$ errors decrease from 45.9\% to 2.7\% for pressure and from 35.7\% to 2.5\% for displacements. 
For snapshot space 2 we do precisely the same experiment with two coarse-grids.
 We present the errors in Table \ref{tab:c1-2} and again see that the errors are also decrease and have roughly the same behavior. 
 In general, we see that the two snapshot choices in this example do not differ greatly and no clear choice arises. In some cases the snapshot space 1 appears to fair better, however, this is not always true. 
Finally, we note that, for solution of nonlinear problem in each time step, the Picard iteration converges after about 3 steps.

\begin{figure}[htb]
\begin{center}
\includegraphics[width=0.45\linewidth]{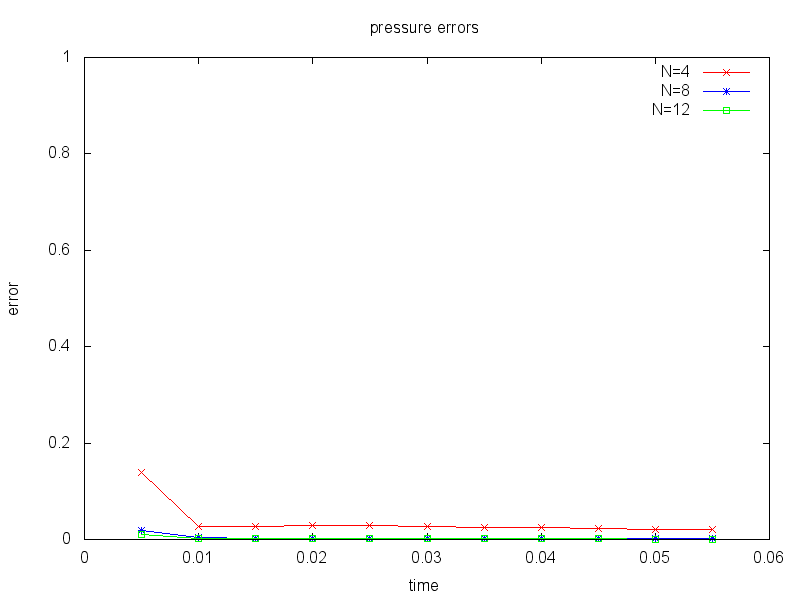} 
\includegraphics[width=0.45\linewidth]{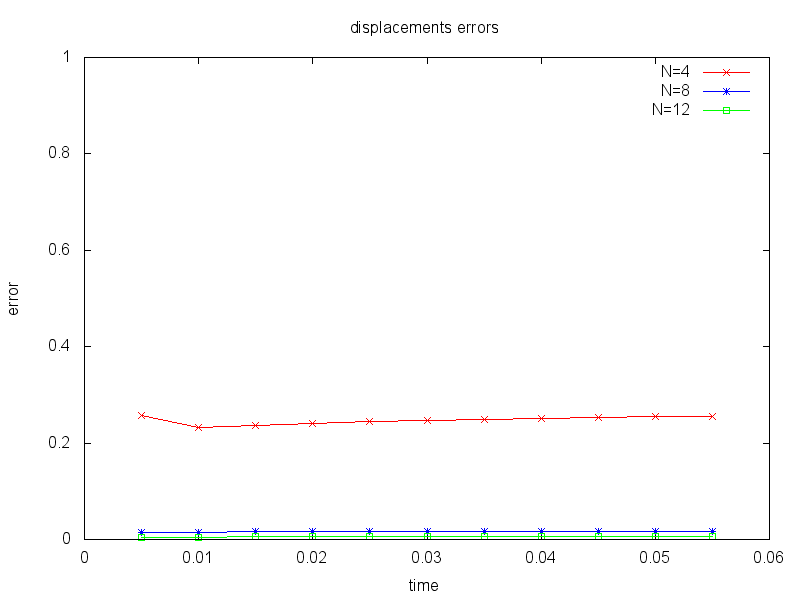} 
\\
\includegraphics[width=0.45\linewidth]{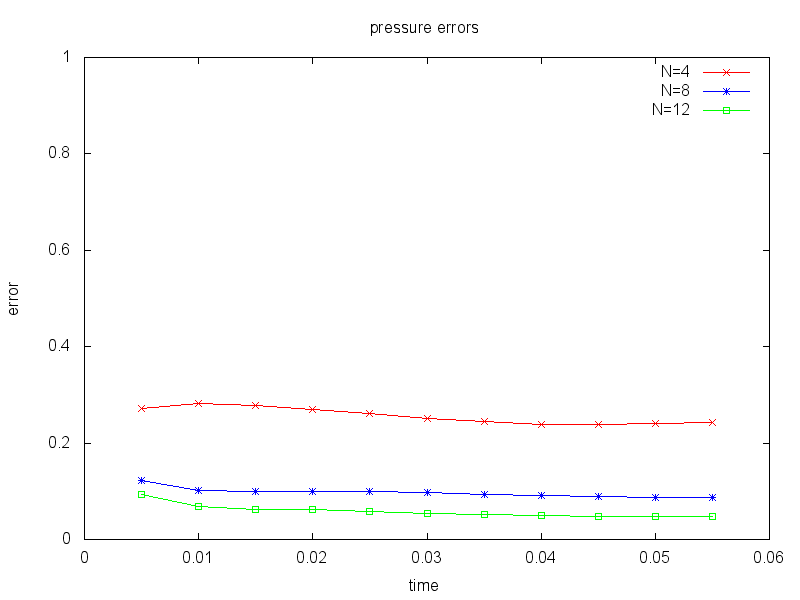} 
\includegraphics[width=0.45\linewidth]{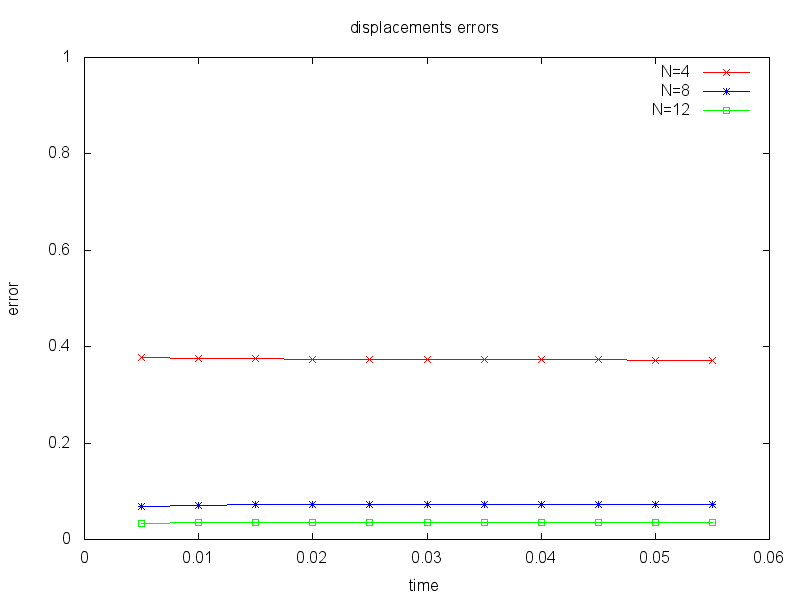} 
\caption{Weighted $L^2$ are on the top and $H^1$ are on the bottom. Errors for pressure are on the left and displacements are on the right  for nonlinear problem on coarse mesh with 36 nodes.}
\label{fig:err-5-1}
\end{center}
\end{figure}

\begin{figure}[htb]
\begin{center}
\includegraphics[width=0.45\linewidth]{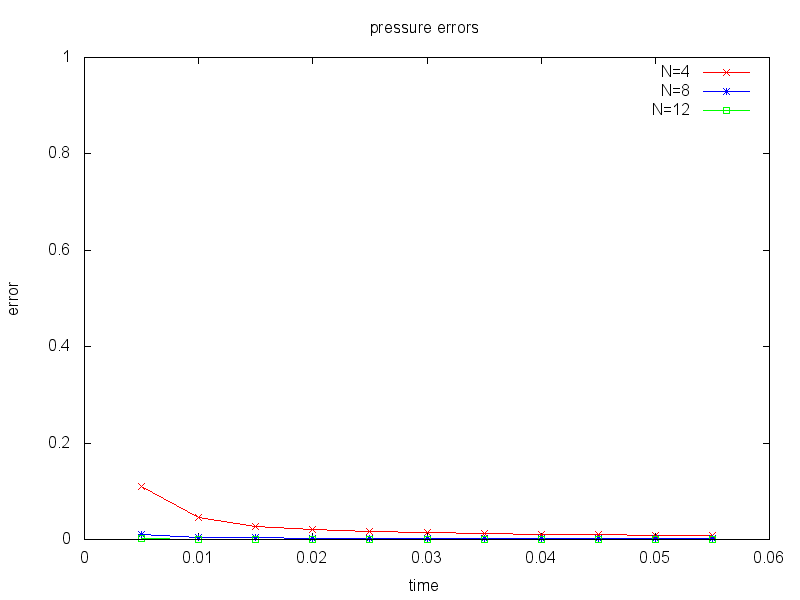} 
\includegraphics[width=0.45\linewidth]{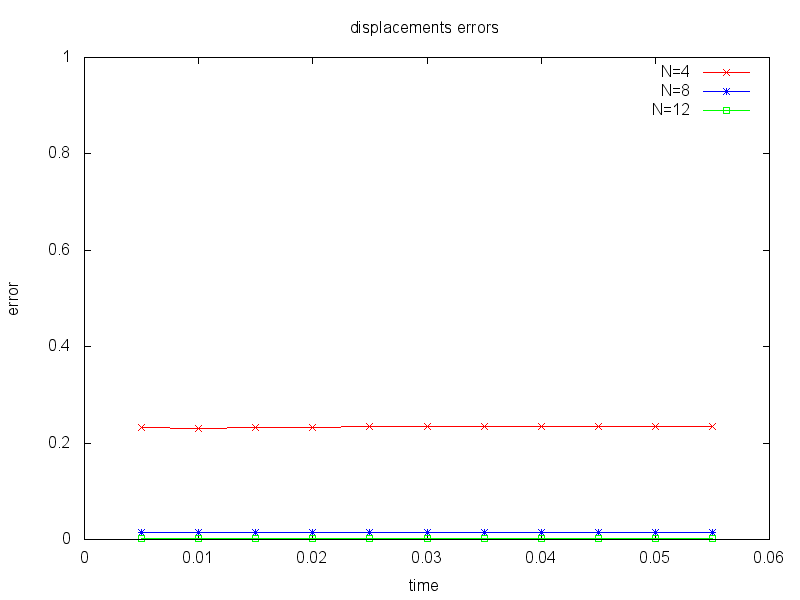} 
\\
\includegraphics[width=0.45\linewidth]{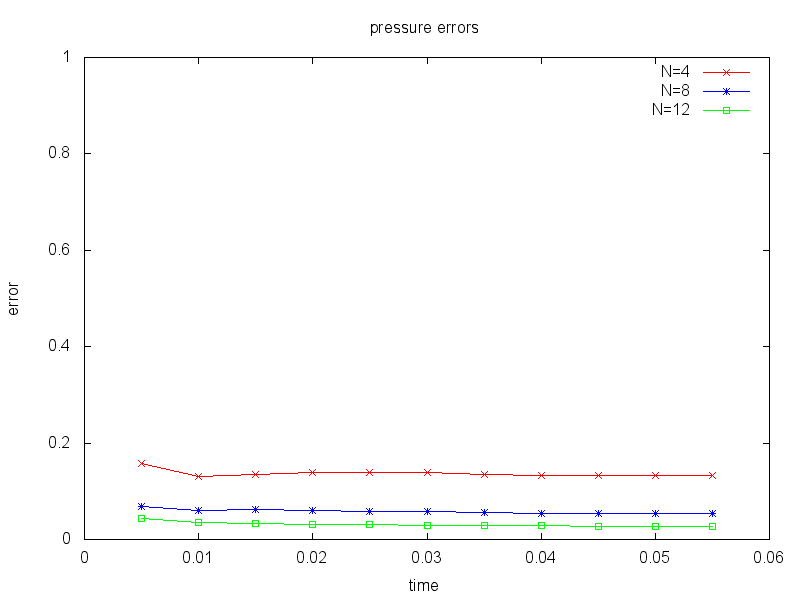} 
\includegraphics[width=0.45\linewidth]{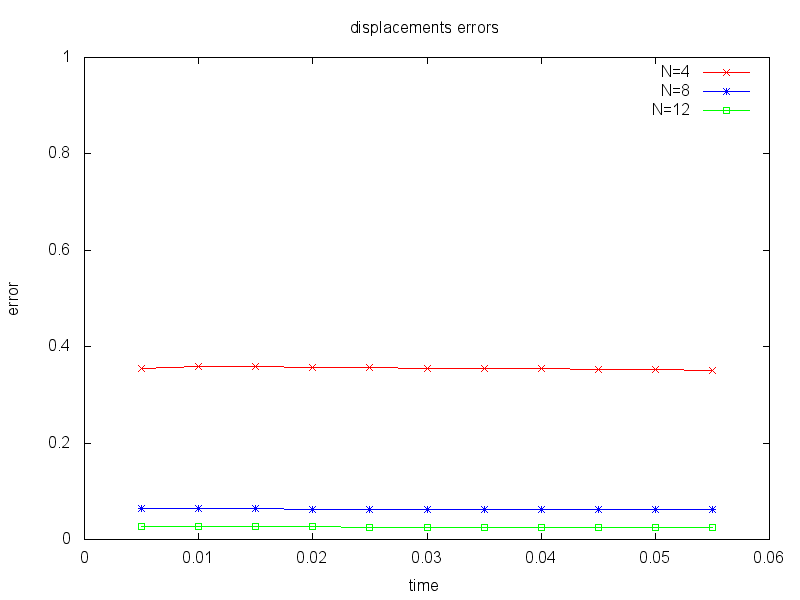} 
\caption{Weighted $L^2$ are on the top and $H^1$ are on the bottom. Errors for pressure are on the left and displacements are on the right  for nonlinear problem on coarse mesh with 121 nodes.}
\label{fig:err-10-1}
\end{center}
\end{figure}

%
To show the stability of the multiscale spaces over time we include time plots. 
We  include plots over time of the error with respect to number of basis functions used. 
To get an idea of the behavior we only present the results for snapshot space 1 for two coarse grids. 
In Figure \ref{fig:err-5-1} and \ref{fig:err-10-1} we show errors over time for $N_{\text{on}} = N_{\text{on}}^{p} = N_{\text{on}}^{u} = 4, 8, 12,$ respectively.
 We observe that errors decrease as we increase the dimension of the offline space as expected and the basis appears to be robust with respect to longer times.

\begin{table}[htp]
\begin{center}
\begin{tabular}[hp]{|c|cc|cc|}
\hline
 & \multicolumn{2}{|c|}{$\varepsilon_p$}
 & \multicolumn{2}{|c|}{$\varepsilon_u$} \\
$N_{\text{on}}^p$  & $L^2$ & $H^1$ & $L^2$ & $H^1$\\
\hline \hline
\multicolumn{5}{|c|}{$N_{\text{on}}^u = 4$}  \\
\hline
2     & 0.081		& 0.609   &  0.304   &  0.380  \\
4     & 0.019		& 0.242   &  0.254   &  0.371  \\
\hline 
\multicolumn{5}{|c|}{$N_{\text{on}}^u = 8$}  \\
\hline
2     & 0.082		& 0.607   &  0.091   &  0.104  \\
4     & 0.021		& 0.241   &  0.023   &  0.074  \\
8     & 0.001		& 0.087   &  0.016   &  0.072  \\
\hline 
\multicolumn{5}{|c|}{$N_{\text{on}}^u= 12$}  \\
\hline
2     & 0.082		& 0.607   &  0.085   &  0.077  \\
4     & 0.021		& 0.241   &  0.013   &  0.037  \\
8     & 0.001		& 0.087   &  0.007   &  0.034  \\
12   & 0.0009	& 0.047   &  0.005   &  0.034  \\
\hline 
\end{tabular}
$\;\;\;$  
\begin{tabular}[hp]{|c|cc|cc|}
\hline
 & \multicolumn{2}{|c|}{$\varepsilon_p$}
 & \multicolumn{2}{|c|}{$\varepsilon_u$} \\
$N_{\text{on}}^p$  & $L^2$ & $H^1$ & $L^2$ & $H^1$\\
\hline \hline
\multicolumn{5}{|c|}{$N_{\text{on}}^u = 4$}  \\
\hline
2     & 0.048 	& 0.459   & 0.264   & 0.357  \\
4     & 0.008 	& 0.132   & 0.235   & 0.351  \\
\hline 
\multicolumn{5}{|c|}{$N_{\text{on}}^u = 8$}  \\
\hline
2     & 0.048 	& 0.457   & 0.063   & 0.079  \\
4     & 0.006 	& 0.130   & 0.022   & 0.063  \\
8     & 0.001 	& 0.053   & 0.015   & 0.062  \\
\hline 
\multicolumn{5}{|c|}{$N_{\text{on}}^u= 12$}  \\
\hline
2     & 0.048 	& 0.457   & 0.052   & 0.051  \\
4     & 0.006 	& 0.130   & 0.009   & 0.026  \\
8     & 0.001 	& 0.053   & 0.002   & 0.025  \\
12   & 0.0002 	& 0.027   & 0.001   & 0.025  \\
\hline 
\end{tabular}
\end{center}
\caption{Numerical results for nonlinear problem using snapshot space 1. Left:  for coarse mesh with 36 nodes. Right:  for coarse mesh with 121 nodes.}
\label{tab:c1-1}
\end{table}


\begin{table}[htp]
\begin{center}
\begin{tabular}[hp]{|c|cc|cc|}
\hline
 & \multicolumn{2}{|c|}{$\varepsilon_p$}
 & \multicolumn{2}{|c|}{$\varepsilon_u$} \\
$N_{\text{on}}^p$  & $L^2$ & $H^1$ & $L^2$ & $H^1$\\
\hline \hline
\multicolumn{5}{|c|}{$N_{\text{on}}^u = 4$}  \\
\hline
2     & 0.063		&  0.551	& 0.302	& 0.380 \\
4     & 0.015		&  0.256	& 0.267	& 0.372 \\
\hline 
\multicolumn{5}{|c|}{$N_{\text{on}}^u = 8$}  \\
\hline
2     & 0.064		& 0.549	& 0.088	& 0.101 \\
4     & 0.013		& 0.255	& 0.036	& 0.078 \\
8     & 0.010		& 0.112	& 0.024	& 0.074 \\
\hline 
\multicolumn{5}{|c|}{$N_{\text{on}}^u= 12$}  \\
\hline
2     & 0.064		& 0.549	& 0.082	& 0.072 \\
4     & 0.013		& 0.255	& 0.028	& 0.040 \\
8     & 0.010		& 0.112	& 0.016	& 0.036 \\
12   & 0.006		& 0.080	& 0.010	& 0.035 \\
\hline 
\end{tabular}
$\;\;\;$  
\begin{tabular}[hp]{|c|cc|cc|}
\hline
 & \multicolumn{2}{|c|}{$\varepsilon_p$}
 & \multicolumn{2}{|c|}{$\varepsilon_u$} \\
$N_{\text{on}}^p$  & $L^2$ & $H^1$ & $L^2$ & $H^1$\\
\hline \hline
\multicolumn{5}{|c|}{$N_{\text{on}}^u = 4$}  \\
\hline
2     & 0.042		& 0.426	& 0.258	& 0.355  \\
4     & 0.008		& 0.145	& 0.235	& 0.351  \\
\hline 
\multicolumn{5}{|c|}{$N_{\text{on}}^u = 8$}  \\
\hline
2     & 0.042		& 0.424	& 0.057	& 0.075  \\
4     & 0.006		& 0.143	& 0.023	& 0.063  \\
8     & 0.001		& 0.078	& 0.015	& 0.062  \\
\hline 
\multicolumn{5}{|c|}{$N_{\text{on}}^u= 12$}  \\
\hline
2     & 0.042		& 0.424	& 0.045	& 0.046  \\
4     & 0.006		& 0.143	& 0.010	& 0.026  \\
8     & 0.001		& 0.078	& 0.002	& 0.025  \\
12   & 0.0001	& 0.039	& 0.001	& 0.025  \\
\hline 
\end{tabular}
\end{center}
\caption{Numerical results for nonlinear problem using snapshot space 2. Left:  for coarse mesh with 36 nodes. Right:  for coarse mesh with 121 nodes.}
\label{tab:c1-2}
\end{table}

\section{Conclusion} 

Modeling and simulation of  a nonlinear poroelastic media  is challenging due the  heterogeneities and the nonlinear dependence on the coefficients.  %
Thus, in this paper we developed a Generalized Multiscale Finite Element Method for a nonlinear poroelastic media. 
%
We gave a general nonlinear poroelasticity model in the framework  of the Biot equations, where we had possibly complex nonlinear dependence on permeability fields and elasticity tensors. 
%
As the Nonlinear GMsFEMs treat nonlinearities as a parameter, we linearize the equations in a time-staggered Picard iteration formulation. 
We then outlined the construction of the multiscale spaces offline and online spaces.
 The algorithm is then implemented on a single geometry with two different cases of permeability fields. The first being the standard linear case  and a second nonlinear relation depending on pressure where a parameter spaces are considered with offline and online spaces. 
 We presented the errors relative to the fine scale solution with varying multiscale basis functions and coarse-grid refinements. Finally, we showed the robustness of the modes for longer time simulations. 


\end{document}